\documentclass[a4paper,11pt]{article}
\usepackage{amsthm, amsmath}

\renewcommand{\baselinestretch}{1}

\makeatletter
  
\def\@seccntformat#1{\csname the#1\endcsname.\quad}
  \makeatother

\makeatletter
  
   \@addtoreset{equation}{section}
\makeatother

\def\<{\left<} \def\>{\right>}
\newtheorem{theorem}{Theorem}

\newtheorem{proposition}[theorem]{Proposition}

\newtheorem{definition}[theorem]{Definition}
\def\proof{\noindent{\it Proof. }}
\def\bea{\begin{eqnarray} }
\def\eea{\end{eqnarray} }
\def\be{\begin{equation} }
\def\ee{\end{equation} }
\def\qed{\ifhmode\unskip\nobreak\fi\ifmmode\ifinner\else\hskip5pt
\fi\fi\hbox{\hskip5 pt \vrule width4 pt height6 pt depth1.5 pt \hskip1pt }}

\begin{document}

\title{Surfaces in Euclidean 3-space whose normal bundles are tangentially biharmonic
\footnote{Archiv der Mathematik {\bf 99} (2012), 281-287. A minor mistake in the proof of Proposition 4 is corrected.}}
\author {Toru Sasahara}


\date{}
\maketitle

\begin{abstract}
{\footnotesize  
A submanifold is said to be tangentially biharmonic if the bitension field of 
the isometric immersion that defines the submanifold  has vanishing tangential component.
The purpose of this paper is to 
prove that a surface in Euclidean $3$-space
 has tangentially biharmonic normal bundle if and
only if it is either minimal, a part of a round sphere, or
a part of a circular cylinder.}
\end{abstract}
\section{Introduction}
A smooth map between Riemannian manifolds is called a biharmonic map if it is a critical point of the bienergy functional.
In \cite{ji2}, Jiang showed that a smooth map is biharmonic if and only if its bitension field vanishes identically.
A submanifold is called a biharmonic submanifold if the isometric immersion that defines the submanifold is  a biharmonic map.
Any minimal submanifold
is biharmonic, but
the converse is not true in general. 
In fact, there exist
many examples of nonminimal biharmonic submanifolds in spheres (see, \cite{bal}). 
Investigating the properties of  biharmonic submanifolds is nowadays becoming a very active field of study. 

On the other hand, Lagrangian submanifolds are the most fundamental objects in symplectic geometry.
The study of them  from the Riemannian geometric point of view was
initiated in 1970's.
In particular, minimal Lagrangian
submanifolds 
have attracted considerable attention from both the geometric and the physical
point of view. 
Some extensions
of such
submanifolds from the view point of 
variational calculus, e.g., Hamiltonian stationary
Lagrangian submanifolds also have been studied widely.

It is known that the normal bundle of a submanifold in Euclidean $n$-space ${\bf E}^n$ can be immersed as a Lagrangian submanifold
in complex Euclidean $n$-space (see, \cite{hl}).   
Harvey and Lawson \cite{hl} showed that a surface in ${\bf E}^3$ 
has minimal normal bundle if and only if it is minimal.
Sakaki \cite{sakaki} generalized their result as follows: a surface in ${\bf E}^3$ has Hamiltonian stationary normal bundle if and only if
it is either minimal, a part of  a round sphere, or a cone with vertex angle $\pi/2$.

In this paper, we obtain a new generalization of Harvey and  Lawson's result mentioned above.
We introduce the notion of  tangentially biharmonicity which is weaker than  
 biharmonicity for submanifolds, and
prove that a surface in ${\bf E}^3$ has tangentially biharmonic normal bundle if and only if it is either minimal, a part of a round sphere, or a circular cylinder. 
This result  also 
provides a partial answer to B.Y.Chen's conjecture, which states that 
any biharmonic submanifold in Euclidean space is minimal.

\section{Preliminaries}
\subsection{Basic formulas}
Let $M^m$ be  
an $m$-dimensional submanifold of a
Euclidean $n$-space ${\bf E}^n$.  
 We denote by $\nabla$ and
$\tilde\nabla$ the Levi-Civita connections on $M^m$ and 
${\bf E}^n$, respectively. The
formulas of Gauss and Weingarten are given respectively by
\begin{equation}
\begin{array}{crl}
& &\tilde \nabla_XY= \nabla_XY+h(X,Y),\\
& &\tilde\nabla_X \xi = -A_{\xi}X+D_X\xi,
\end{array}
\end{equation}
 for tangent vector fields $X$, $Y$ and a normal vector field $\xi$, 
where $h,A$ and $D$ are the second fundamental
form, the shape operator and the normal
connection.
The shape operator and the second fundamental form are related by
\be
\langle h(X, Y), \xi \rangle=\langle A_\xi X, Y\rangle.
\ee
The mean curvature vector field $H$ is defined by 
$H=(1/m){\rm trace}\hskip3pt h$.
If it vanishes identically, then $M^m$ is called a {\it minimal submanifold}.

The equation of Codazzi is given by
\be
({\bar\nabla}_{X}h)(Y,Z)=
({\bar\nabla}_{Y}h)(X,Z),
\ee
where $X,Y,Z,W$ 
 are vectors tangent  to
$M^m$,  and
$\bar\nabla h$ is defined by
\be ({\bar\nabla}_{X}h)(Y,Z)= D_X h(Y,Z) - h(\nabla_X
Y,Z) - h(Y,\nabla_X Z).\nonumber
\ee
In this paper, we do not use the equations  of Gauss and Ricci. Therefore, we omit them.

\subsection{Tangentially biharmonic submanifolds}
Let $f:(M^m, g) \rightarrow N$ be a smooth map between  two Riemannian manifolds.
The {\it tension field}
$\tau(f)$ of $f$ is a section
of the induced vector bundle $f^{*}TN$
defined by
\be
\tau(f):=\mathrm{trace}_g(\nabla^{f} df)
\ee
where $\nabla^f$ denotes the induced connection of $M$.
A smooth map $f$ is called a {\it harmonic map} if $\tau(f)$ vanishes identically.
In the case where $f$ is an isometric immersion,  we have 
\be\tau(f)=mH.\ee

The bienergy $E_2(f)$ of $f$ over compact domain $\Omega\subset M^m$ is defined by 
\be
E_2(f)=\int_{\Omega}|\tau(f)|^2dv_g,\nonumber
\ee
where $dv_g$ is the volume form of $M^m$ (see, \cite{es2}).
Then $E_2$ provides a measure for the extent to which $f$ fails to be harmonic. 
If $f$ is a critical point of $E_2$
over every compact domain, then $f$ is called a
{\it biharmonic map} (or {\it $2$-harmonic map}).
In \cite{ji2}, Jiang proved that $f$ is biharmonic if and only if its bitension field defined by
\be\tau_2(f):=-\Delta_f\tau(f)+{\rm trace}_gR^{N}(\tau(f),df)df\ee
vanishes identically, where
$\Delta_f=-{\rm trace}_g(\nabla^{f}\nabla^{f}-\nabla^{f}_{\nabla})$ and 
$R^{N}$ is 
the curvature tensor  of $N$, which is given by $$R^{N}(X, Y)Z=[\nabla^N_X, \nabla^N_Y]Z-\nabla^N_{[X, Y]}Z$$ for the Levi-Civita connection $\nabla^N$ of $N$.

If $f$ is an isometric immersion, then $M^m$ or $f(M^m)$ is called a {\it biharmonic submanifold} in $N$. 
It follows from (2.5) and (2.6) that any minimal submanifold is biharmonic. However, the converse is not true in general. In fact, there exists many nonminimal biharmonic submanifolds in spheres (see \cite{bal}).
On the other hand, B.Y.Chen \cite{chen} proposed the conjecture: any biharmonic submanifold in  Euclidean space is minimal. Several partial answers have been obtained, but the conjecture is still open in general. 

We introduce the notion of  tangentially biharmonicity which is weaker than 
 biharmonicity for submanifolds as follows:
\begin{definition}
A submanifold  is said to be {\it tangentially biharmonic}
if it satisfies  
\be
\{\tau_2(f)\}^{\top}=0,
\ee where $\{\cdot\}^\top$ denotes the tangential part of $\{\cdot\}$. 
\end{definition}

\subsection{The tension field of normal bundle of surfaces in ${\bf E}^3$}

Let $x:M^2\rightarrow {\bf E}^3$ be an isometric immersion from a Riemannian 2-manifold into Euclidean 3-space.
The normal bundle $T^{\perp}M^2$ of $M^2$ is naturally immersed in ${\bf E}^3\times{\bf E}^3={\bf E}^6$ by the immersion 
$f(\xi_x):=(x, \xi_x)$, which is  expressed as
\be
f(x, t)=(x, tN)
\ee
for the unit normal vector field $N$ along $x$.

We equip $T^{\perp}M^2$ with the metric induced by $f$. 
If we define the complex structure $J$ on ${\bf C}^3={\bf E}^3\times{\bf E}^3$
by $J(X, Y):=(-Y, X)$, then 
$T^{\perp}M^2$ is a Lagrangian submanifold in ${\bf C}^3$, 
that is, it satisfies $\langle JX, Y\rangle=0$ for any two vectors $X$ and $Y$ tangent to
$T^{\perp}M^2$  (see [4, III.3.C]). 


In \cite{sakaki}, the tension field  $\tau(f)$ of $f$ is represented by  
a local coordinate system such that the coordinate curves
are lines of curvature on $M^2$. 
In this section, we  reexpress  $\tau(f)$ by an orthonormal frame of principal directions on $M^2$  for later use.

Let $e_1$ and $e_2$ be principal directions with 
the corresponding principal curvatures $a$ and $b$, respectively,
 that is,
\begin{equation}
\begin{array}{crl}
A_Ne_1=ae_1,\quad
A_Ne_2=be_2.
\end{array}
\end{equation}
Put 
\begin{equation}
\begin{array}{crl}
&&\tilde e_1=(1+t^2a^2)^{-\frac{1}{2}}e_1, \\
&&\tilde e_2=(1+t^2b^2)^{-\frac{1}{2}}e_2, \\
&&\tilde e_3=\frac{\partial}{\partial t}.
\end{array}
\end{equation}
Then, it follows from (2.1), (2.8) and (2.9) that
\begin{equation}
\begin{array}{crl}
&&f_{*}(\tilde e_1)=(1+t^2a^2)^{-\frac{1}{2}}(e_1, -tae_1),\\
&&f_{*}(\tilde e_2)=(1+t^2b^2)^{-\frac{1}{2}}(e_2, -tbe_2),\\
&&f_{*}(\tilde e_3)=(0, N).
\end{array}
\end{equation}
We see that $\{\tilde e_1, \tilde e_2, \tilde e_3\}$ is an orthonormal frame on $T^{\perp}M^2$ .
We set
\begin{equation}
\begin{array}{crl}
&&e_4:=Jf_{*}(\tilde e_1)=(1+t^2a^2)^{-\frac{1}{2}}(tae_1, e_1),\\
&&e_5:=Jf_{*}(\tilde e_2)=(1+t^2b^2)^{-\frac{1}{2}}(tbe_2, e_2),\\
&&e_6:=Jf_{*}(\tilde e_3)=(-N, 0).
\end{array}
\end{equation}
Then $\{e_4, e_5, e_6\}$ is a normal orthonormal frame along $f$. 

Put $\<\nabla_{e_i}e_j, e_k\>=\omega_j^k(e_j)$ for $i, j ,k\in\{1, 2\}$. Note that $\omega_1^2(e_k)=
-\omega_2^1(e_k)$.
The equation (2.3) of Codazzi yields
\begin{equation}
\begin{array}{crl}
e_1b=(a-b)\omega_1^2(e_2),\\
e_2a=(b-a)\omega_2^1(e_1).
\end{array}
\end{equation}

We put
$h^{\alpha}_{ij}=\langle{\tilde e_i}(f_{*}(\tilde e_j)), e_{\alpha}\rangle$, $1\leq i, j\leq 3$, $4\leq\alpha\leq 6$.
Then, (2.1) and (2.5) imply that $\tau(f)$  is given by 
\be
\tau(f)=\sum_{\alpha=4}^{6}\sum_{i=1}^3h_{ii}^{\alpha}e_\alpha.
\ee
It follows from (2.10)-(2.13) that
\begin{equation}
\begin{array}{crl}
&&h^4_{11}=-t(1+t^2a^2)^{-\frac{3}{2}}e_1a,\\
&&h^4_{22}=-t(1+t^2a^2)^{-\frac{1}{2}}(1+t^2b^2)^{-1}(b-a)\omega_2^1(e_2),\\
&&\hskip19pt=-t(1+t^2a^2)^{-\frac{1}{2}}(1+t^2b^2)^{-1}e_1b,\\
&&h^5_{11}=-t(1+t^2a^2)^{-1}(1+t^2b^2)^{-\frac{1}{2}}(a-b)\omega_1^2(e_1),\\
&&\hskip19pt=-t(1+t^2a^2)^{-1}(1+t^2b^2)^{-\frac{1}{2}}e_2a,\\
&&h^5_{22}=-t(1+t^2b^2)^{-\frac{3}{2}}e_2b,\\
&&h^{6}_{11}=-a(1+t^2a^2)^{-1},\\
&&h^{6}_{22}=-b(1+t^2b^2)^{-1},\\
&&h^{4}_{33}=h^5_{33}=h^6_{33}=0.
\end{array}
\end{equation}
Set
\begin{equation}
\begin{array}{crl}
&&P=(1+t^2a^2)^{-2}e_1a+(1+t^2a^2)^{-1}(1+t^2b^2)^{-1}e_1b,\\
&&Q=(1+t^2a^2)^{-1}(1+t^2b^2)^{-1}e_2a+(1+t^2b^2)^{-2}e_2b,\\
&&R=a(1+t^2a^2)^{-1}+b(1+t^2b^2)^{-1}.
\end{array}
\end{equation}
Then by (2.12), (2.14) and (2.15) we find (cf. \cite{sakaki})
\be
\tau(f)=-(Pt^2ae_1+Qt^2be_2-RN, Pte_1+Qte_2).
\ee

Using (2.16) and (2.17) we obtain the following (cf. [4, III. Th.3.11, Pop.2.17]):
\begin{proposition}
 A surface $M^2$ in ${\bf E}^3$ is minimal if and only if 
$T^{\perp}M^2$ is a minimal submanifold of ${\bf E}^6$.
\end{proposition}

\section{Main results}

In view of Proposition 2, it is natural to investigate the  biharmonicity of $T^{\perp}M^2$ in ${\bf E}^6$.
The main theorem of this paper is the following:
\begin{theorem}
A surface in Euclidean $3$-space has tangentially biharmonic normal bundle if and only if it is either minimal, a part of 
a round sphere, or
a part of a circular cylinder.
\end{theorem}

\proof
If $M^2$ is minimal, then 
it follows from Proposition 2 that $T^{\perp}M^2$ is tangentially biharmonic.
Assume that the mean curvature vector field of $M^2$ vanishes nowhere
 and  $T^{\perp}M^2$ is tangentially biharmonic.

Let $(X, Y)$ be  a section of $f^{*}T({\bf E}^3\times{\bf E}^3)$
and let $\nabla$ be a differential operator acting on $X$ and $Y$ which is given by 
\bea
\Delta
&=&-(1+t^2a^2)^{-1}(\tilde\nabla_{e_1}\tilde\nabla_{e_1}-\tilde\nabla_{\nabla_{e_1}e_1})\nonumber\\
&&-(1+t^2b^2)^{-1}(\tilde\nabla_{e_2}\tilde\nabla_{e_2}-\tilde\nabla_{\nabla_{e_2}e_2})-\frac{\partial^2}{\partial t^2}
\eea
for the Levi-Civita connection $\tilde\nabla$ of ${\bf E}^3$. 
Then, we have
\be\Delta_f(X, Y)=(\Delta X, \Delta Y).\ee

First, using (2.6), (2.17) and (3.2), we calculate $\tau_2(f)|_{t=0}$.
By (2.16) and (3.1) we get
\bea
&&\Delta(Pt^2ae_1+Qt^2be_2)|_{t=0}=-2ae_1(a+b)e_1-2be_2(a+b)e_2,\\
&&\Delta(Pte_1+Qte_2)|_{t=0}=0.
\eea
Moreover, it follows from  (2.1), (2.2), (2.9), (2.13) and (3.1) that
\bea
\Delta N|_{t=0}&=&-\tilde\nabla_{e_1}\tilde\nabla_{e_1}N+\tilde\nabla_{\nabla_{e_1}{e_1}}N
-\tilde\nabla_{e_2}\tilde\nabla_{e_2}N+\tilde\nabla_{\nabla_{e_2}{e_2}}N\nonumber\\
&=&\tilde\nabla_{e_1}(ae_1)-\omega_1^2(e_1)be_2+\tilde\nabla_{e_2}(be_2)-\omega_2^1(e_2)ae_1\nonumber\\
&=&\{e_1a+(b-a)\omega_2^1(e_2)\}e_1+\{e_2b+(a-b)\omega_1^2(e_1)\}e_2\nonumber\\
&&\quad+ah(e_1, e_1)+bh(e_2, e_2)\nonumber\\
&=&e_1(a+b)e_1+e_2(a+b)e_2+(a^2+b^2)N.
\eea
By (2.16), (3.1) and (3.5), we obtain
\bea
\Delta(RN)|_{t=0}&=&\Delta_M(a+b)N-2e_1(a+b)\tilde\nabla_{e_1}N-2e_2(a+b)\tilde\nabla_{e_2}N+(a+b)\Delta N\nonumber\\
&&-R_{tt}|_{t=0}N\nonumber\\
&=&\Delta_M(a+b)N+(3a+b)e_1(a+b)e_1+(a+3b)e_2(a+b)e_2\nonumber\\
&&+(a+b)(3a^2-2ab+3b^2)N,
\eea
where $\Delta_M\phi=-\sum_{i=1}^2\{e_i(e_i\phi)-(\nabla_{e_i}e_i)\phi\}$ for a smooth function $\phi$ of $M^2$.
Consequently, it follows from (2.6), (2.17), (3.2), (3.3), (3.4) and (3.6) that 
\bea
\tau_2(f)|_{t=0}&=&-(\Delta_M(a+b)N+(5a+b)e_1(a+b)e_1+(a+5b)e_2(a+b)e_2\nonumber\\
&&+(a+b)(3a^2-2ab+3b^2)N, 0).
\eea
By the assumption, we have
$$\<\tau_2(f), f_{*}(e_1)\>|_{t=0}=\<\tau_2(f), f_{*}(e_2)\>|_{t=0}=0.$$
Therefore, we derive from (2.11) and (3.7) that
\begin{equation}
\begin{array}{crl}
&&(5a+b)e_1(a+b)=0,\\
&&(a+5b)e_2(a+b)=0.
\end{array}
\end{equation}

Next, we compute $\<\tau_2(f), \tilde e_3\>$ at each point on $T^{\perp}M^2$.
By (2.11) and (2.17) we get
\bea
&& \<\tau_2(f), \tilde e_3\>
=\<\Delta(Pte_1+Qte_2),N\>\nonumber\\
&=&-t(1+t^2a^2)^{-1}\{2(e_1P)a-\omega_1^2(e_1)Qb+P\langle\tilde\nabla_{e_1}\tilde\nabla_{e_1}e_1, N\rangle+Q\langle\tilde\nabla_{e_1}\tilde\nabla_{e_1}e_2, N\rangle\}\nonumber\\
&&-t(1+t^2b^2)^{-1}\{2(e_2Q)b-\omega_2^1(e_2)Pa+P\langle\tilde\nabla_{e_2}\tilde\nabla_{e_2}e_1, N\rangle
+Q\langle\tilde\nabla_{e_2}\tilde\nabla_{e_2}e_2, N\rangle\}.
\eea
Substituting (2.16) into (3.9), we find that 
\bea
&&(1+t^2a^2)^4(1+t^2b^2)^4\<\tau_2(f), \tilde e_3\>\nonumber\\
&=&8(1+t^2b^2)^4t^3a^2(e_1a)^2+
8(1+t^2a^2)^4t^3b^2(e_2b)^2+(1+t^2a^2)(1+t^2b^2)B
\eea
for some function $B$ on $T^{\perp}M^2$, which is polynomial with respect to $t$.

The assumption yields that (3.10) is identically zero, also as polynomial for $t\in{\bf C}$. 
If $ab\ne 0$, then
substituting $t=\sqrt{-1}/a$ and $t
=\sqrt{-1}/b$ into (3.10) shows 
\be
(a-b)^2e_1a=(a-b)^2e_2b=0.
\ee
Combining (3.8) with (3.11), we state that $a$ and $b$ are nonzero constant and hence $M^2$ is a part of a round sphere. 
If $ab=0$, then we may assume that $a\ne 0$ and $b=0$. By (3.8), we conclude that $a$ is constant and hence $M^2$
is a part of a circular cylinder.

Conversely, if $M^2$ is either a part of a round sphere or a part of a circular cylinder, then $a$ and $b$ are constant, hence $P=Q=0$ from (2.16). By using (2.17), (3.1) and (3.2) we can see that $T^{\perp}M^2$ is nonminimal  and tangentially biharmonic.\qed\vspace{1.5ex}
\renewcommand{\baselinestretch}{1}

It follows from (3.7) that
 the normal bundles of a round sphere and  a  circular cylinder cannot be biharmonic.
  Therefore, by Proposotion 2 and Theorem 3
we obtain a partial answer to B. Y. Chen's conjecture as follows:
\begin{proposition}
The normal bundle of a surface in ${\bf E}^3$ is a biharmonic submanifold in ${\bf E}^6$
 if and only if it is minimal.
\end{proposition}


\medskip
{\renewcommand{\baselinestretch}{1}
\noindent General Education and Research Center

\noindent Hachinohe Institute
of Technology

\noindent Hachinohe 031-8501, JAPAN

\noindent E-mail address:  {\tt sasahara@hi-tech.ac.jp}

\end{document}